\documentclass[12pt]{article}

\usepackage{lineno}
\nolinenumbers

\usepackage[utf8]{inputenc}
\usepackage[english]{babel}
\usepackage[table]{xcolor}
\usepackage[fleqn]{amsmath}
\usepackage{nccmath}
\usepackage{booktabs}
\usepackage{caption}
\usepackage{float}
\usepackage{sidecap}
\usepackage{csquotes}
\usepackage{scalerel}

\usepackage[backend=biber,style=alphabetic,sorting=ynt]{biblatex}

\usepackage{answers}
\usepackage{setspace}
\usepackage{graphicx}
\usepackage{enumitem}
\usepackage{multicol}
\usepackage{mathrsfs}
\usepackage[margin=1in]{geometry} 
\usepackage{amsmath,amsthm,amssymb}
\usepackage{mathtools}
\usepackage{adjustbox}

\newcommand{\R}{\mathbb{R}}

\newcommand{\tstirl}[2]{ \begin{bsmallmatrix} #1\\ #2 \end{bsmallmatrix} }
\newcommand\asseq{\stackrel{\mathclap{\normalfont\mbox{!}}}{=}}
\newcommand\eqVerm{\stackrel{\mathclap{\normalfont\mbox{(8c)}}}{=}}

\DeclareMathOperator{\Li}{Li}

\usepackage{listings}
\usepackage{xcolor}

\definecolor{codegreen}{rgb}{0,0.6,0}
\definecolor{codegray}{rgb}{0.5,0.5,0.5}
\definecolor{codepurple}{rgb}{0.58,0,0.82}
\definecolor{backcolour}{rgb}{0.95,0.95,0.92}

\lstdefinestyle{mystyle}{
    backgroundcolor=\color{backcolour},   
    commentstyle=\color{codegreen},
    keywordstyle=\color{magenta},
    numberstyle=\tiny\color{codegray},
    stringstyle=\color{codepurple},
    basicstyle=\ttfamily\footnotesize,
    breakatwhitespace=false,         
    breaklines=true,                 
    captionpos=b,                    
    keepspaces=true,                 
    numbers=left,                    
    numbersep=5pt,                  
    showspaces=false,                
    showstringspaces=false,
    showtabs=false,                  
    tabsize=2
}

\begin{document}

\title{On convolution powers of $1/x$}
\author{Andreas B.G. Blobel\\ 
andreas.blobel@kabelmail.de} 

\maketitle

\begin{abstract}
\noindent Convolution powers of $1/x$ are transformed into functions $f_n$, which satisfy a simple recurrence relation. Solutions are characterized and analyzed.
\end{abstract}
\vspace{10mm}

\setcounter{secnumdepth}{0} 
\setcounter{tocdepth}{4} 

\tableofcontents
\listoftables

\subsection{0. Overview}
In section 1 it is shown, that convolution powers of $1/x$ may be transformed into functions $f_n$, as given in \eqref{eq_237_ansatz_rueck}/\eqref{eq_237_ansatz}.
It is proven by induction, that the $f_n$ solve the recurrence relation
\eqref{eq_f01_repeat}/\eqref{eq_fn0_repeat}/\eqref{eq_Vermutung_237A}
and satisfy the reflection property \eqref{eq_248A_reflection}.

\vspace{5mm}
\noindent
Section 2 contains some definitions used in sections 3 and 4.

\vspace{5mm}
\noindent
In section 3 it is shown, that the functions $f_n$
can be expressed as linear combinations \eqref{eq_A149} of functions $J^{n-k}[1]$, with coefficients $\beta_k$ being determined by the recurrence relation \eqref{eq_beta0}/\eqref{eq_betaRec_150A}.

\vspace{5mm}
\noindent
Section 4 is dedicated to the further analysis of the functions $J^n[1]$,
$n$-fold powers of the $J$ operator \eqref{eqJ} applied to the $1$-function.
They can be decomposed \eqref{eqJm_154B} into expressions of the type $\ Q_j(x) \,\cdot\, \ln(x)^{n-j}\ $, with functions $Q_j$ being generated by the recurrence relation \eqref{eq_Qini}/\eqref{eq_Qrec_160A}.
The functions $Q_j$, in turn, are power series \eqref{eq_Q} in $1/x$, whose coefficients are determined by \eqref{eq_q251B_a}/\eqref{eq_q251B_b}.
Here, lower triangular integer matrices $A^s_{m,\,j}$ become involved, which may be computed from \eqref{eq_q236A_a}/\eqref{eq_q236A_b}. Example data are listed in Table \ref{list:A}.

\subsection{1. Convolution powers of $1/x$}
\noindent
Let $\lambda$ and $a$ be real parameters fulfilling the condition
\begin{ceqn}
\begin{align}
\lambda + a \;>\; 0
\label{eq_237_gt0}
\end{align}
\end{ceqn}
Define real functions $\varphi_{\lambda,a}: \R \rightarrow \R$
\begin{subequations}
\begin{align}
\varphi_{\lambda,a}(x) &=\; 0              & x < \lambda
\label{eq_237_varphi0}
\\[8pt]
\varphi_{\lambda,a}(x) &=\; \frac{1}{x+a}  & x \ge \lambda
\label{eq_237_varphi}
\end{align}
\end{subequations}
Convolution powers \cite{ConvolutionPower} of $\varphi_{\lambda,a}$ are defined by the recurrence relation
\begin{subequations}
\begin{align}
\varphi_{\lambda,a}^{*1} &\;=\; \varphi_{\lambda,a}
\label{eq_237_convolution_init}\\[8pt]
\varphi_{\lambda,a}^{*(n+1)}(x) &\;=\;
\int_{-\infty}^\infty \varphi_{\lambda,a}(x-t) \;\cdot\; \varphi_{\lambda,a}^{*n}(t) \; dt
\hspace{68pt}  x \in \R  \hspace{49pt}  n \ge 1
\label{eq_237_convolution}
\end{align}
\end{subequations}
\newline\noindent
Cut-off condition \eqref{eq_237_varphi0} translates into the general property
\begin{align}
\varphi_{\lambda,a}^{*n}(x) \;&=\; 0   & x < n \cdot \lambda \hspace{29pt}  n \ge 1
\label{eq_27B}
\end{align}
Regarding the non-zero part of $\varphi_{\lambda,a}^{*n}(x)$, consider the transformation
\begin{subequations}
\begin{align}
\varphi_{\lambda,a}^{*n}(x) \;&=\;\;\, \frac{n!}{x + n \cdot a}
\cdot\; f_{n-1}\Big(\,\frac{x-n \cdot \lambda}{\lambda+a}\Big)
 & x \ge n \cdot \lambda \hspace{29pt}  n \ge 1
\label{eq_237_ansatz_rueck}
\\[8pt]
f_{n-1}(y) \;&=\;  \frac{y+n}{n!}\, \cdot (\lambda + a)
\cdot  \varphi_{\lambda,a}^{*n}\big(\, (\lambda + a)\cdot y + n \cdot \lambda\big)
 & y \ge 0  \hspace{46pt}  n \ge 1
\label{eq_237_ansatz}
\end{align}
\end{subequations}
The functions $f_n$ are defined on non-negative real numbers, and their index $n$ starts from zero.

\vspace{1cm}
\noindent
Rewriting \eqref{eq_237_convolution_init}/\eqref{eq_237_convolution} in terms of $f_n$ yields
\begin{subequations}
\begin{align}
f_o \;&=\; 1
\label{eq_f01}\\[8pt]
(n+1) \cdot f_n(y) \;&=\; (y+n+1) \cdot \int_o^y \frac{1}{y-s+1} \cdot \frac{1}{s+n} \cdot f_{n-1}(s) \; ds
\nonumber\\[8pt]
&=\; \int_o^y \frac{1}{y-s+1} \cdot f_{n-1}(s) \; ds 
\;\;+\; \int_o^y \frac{1}{s+n} \cdot f_{n-1}(s) \; ds
\nonumber\\[12pt]
& \hspace{19pt}  y \ge 0 \hspace{30pt}  n \ge 1
\label{eq_248_fRec}
    \end{align}
\end{subequations}
This makes clear, that through the transformation \eqref{eq_237_ansatz_rueck}/\eqref{eq_237_ansatz}, the parameters $\lambda$ and $a$ have been eliminated.
\eqref{eq_248_fRec} implies
\begin{align}
f_n(0) \;&=\; 0
\hspace{98pt}  n \ge 1
\label{eq_fn0}
\end{align}

\vspace{1cm}
\noindent
Consider the set of conditions
\begin{subequations}
\begin{align}
f_o \;&=\; 1
\label{eq_f01_repeat}\\[8pt]
f_n(0) \;&=\; 0 & \hspace{29pt}  n \ge 1
\label{eq_fn0_repeat}\\[8pt]
f^\prime_n(y) \;&\asseq\; \frac{1}{y+n} \;\cdot\;  f_{n-1}(y)  &  y \ge 0 \hspace{29pt}  n \ge 1
\label{eq_Vermutung_237A}
\end{align}
\end{subequations}
\eqref{eq_f01_repeat} and \eqref{eq_fn0_repeat} just repeat \eqref{eq_f01} and  \eqref{eq_fn0}.
In the recurrence relation \eqref{eq_Vermutung_237A}, $f^\prime_n$ denotes the \textit{derivative} of $f_n$ 
with respect to its argument.

\vspace{5mm}
\noindent
Assertion \eqref{eq_Vermutung_237A} will be proved by induction.
Replacing $n=1$ in \eqref{eq_248_fRec}, while observing \eqref{eq_f01}, gives
\begin{align}
\label{eq_f1}
2 \cdot f_1(y) \;&=\; \int_o^y \frac{1}{y-s+1} \; ds 
\;\;+\; \int_o^y \frac{1}{s+1} \; ds
\nonumber\\[8pt]
&\;=\; \bigg[-\ln(y-s+1) \;+\; \ln(s+1)\bigg]^{s=y}_{s=o}
\nonumber\\[8pt]
&\;=\; 2 \cdot \ln(y+1) & y \ge 0
\end{align}
which verifies the induction hypothesis \eqref{eq_Vermutung_237A} for $n=1$.

\vspace{5mm}
\noindent
When deriving \eqref{eq_248_fRec} with respect to $y$, while applying the \textit{Leibniz} integral rule \cite{LeibnizIntegralRule}, one gets
\begin{align}
(n+1) \cdot f^\prime_n(y) \;&=\; \Big( 1 + \frac{1}{y+n} \Big) \cdot f_{n-1}(y)
\nonumber\\[8pt]
&\;-\; \int_o^y \frac{1}{(y-s+1)^2} \cdot f_{n-1}(s) \;ds
\hspace{59pt}  y \ge 0 \hspace{29pt}  n \ge 1
\label{eq_249_2}
\end{align}
Substituting $n \rightarrow n+1$, and applying partial integration to the second term with respect to $s$ yields
\begin{align}
(n+2) \cdot f^\prime_{n+1}(y) \;&=\; \frac{1}{y+n+1} \cdot f_{n}(y)
 \;+\;  \frac{1}{y+1} \cdot f_{n}(0)
\nonumber\\[8pt]
&\;+\; \int_o^y \frac{1}{y-s+1} \cdot f^\prime_{n}(s) \;ds
\hspace{75pt}  y \ge 0 \hspace{29pt}  n \ge 0
\label{eq_249_3}
\end{align}
On the other hand, using the induction hypothesis \eqref{eq_Vermutung_237A} and \eqref{eq_248_fRec}, we have
\begin{align}
\int_o^y \frac{1}{y-s+1} \cdot f^\prime_{n}(s) \;ds
\;&=\; \int_o^y \frac{1}{y-s+1} \cdot \frac{1}{s+n} \cdot f_{n-1}(s) \;ds
\nonumber\\[8pt]
\;&=\;  \frac{n+1}{y+n+1} \cdot f_n(y)
&  y \ge 0 \hspace{29pt}  n \ge 1
\label{eq_249_4}
\end{align}
Finally, inserting \eqref{eq_249_4} into  \eqref{eq_249_3} gives
\begin{align}
f^\prime_{n+1}(y) \;&=\; \frac{1}{n+2} \cdot \bigg( \frac{n+2}{y+n+1} \cdot f_{n}(y)
\;+\;  \frac{1}{y+1} \cdot f_{n}(0) \bigg)
&  y \ge 0 \hspace{29pt}  n \ge 1
\label{eq_249unten}
\end{align}
from which the induction step follows, if one observes \eqref{eq_fn0_repeat}. $\hspace{1cm}\square$

\vspace{8mm}
\noindent
The (non-recursive) \textit{reflection} property
\begin{align}
     \int_o^y \frac{1}{y-s+1} \cdot f_{n}(s) \; ds \;=\; (n+1) \cdot
     \int_o^y \frac{1}{s+n+1} \cdot f_{n}(s) \; ds  \hspace{30pt}  y \ge 0 \hspace{22pt}  n \ge 0
    \label{eq_248A_reflection}
\end{align}
follows as a corollary, if one replaces
\begin{ceqn}
\begin{equation*}
f_n(y) \;=\; \int_o^y f^\prime_{n}(s) \; ds \;\ \eqVerm\ \; \int_o^y \frac{1}{s+n} \;\cdot\;  f_{n-1}(s) \; ds
\end{equation*}
\end{ceqn}
in the left-hand side of \eqref{eq_248_fRec}.

\subsection{2. Some definitions}
Given a smooth, real function $f$, declare the `harmonic' integration operator $H$ as the antiderivative \cite{Antiderivative}
\begin{align} \label{eq_H}
H[f] &:= \int \frac{f(x)}{x} \; dx
\end{align}
In particular, if $g$ denotes a power series in $1/x$ with radius of convergence $C > 0$
\begin{align} \label{eq_g}
g(x) &= \sum_{k = 0}^\infty a_k \cdot x^{-k}    & x > 1/C > 0
\end{align}
one gets
\begin{align} \label{eq_H1g}
H[g](x) &= a_o \cdot \ln(x) \; - \; \sum_{k = 1}^\infty \, \frac{a_k}{k} \cdot x^{-k}    & x > 1/C
\end{align}
Here, $\ln$ denotes the natural logarithm.
Clearly, if $g(x)$ converges for $x > 1/C$, so does the second term in \eqref{eq_H1g}.
\newline
More generally, applying $H$ multiple times, we get
\begin{align} \label{eq_Hmg}
H^m[g](x) &= a_o \cdot \frac{\ln(x)^m}{m!} \; + \; (-1)^m \; \cdot \; \sum_{k = 1}^\infty \, \frac{a_k}{k^m} \cdot x^{-k}    & m \ge 0 \hspace{29pt} x > 1/C
\end{align}
Let $\nabla$ denote the \textit{backward difference} operator \cite{FiniteDifference}
\begin{align}
\nabla[f](x) &:= f(x) - f(x-1)   \label{eq_N}
\end{align}
Applying $\nabla$ to \eqref{eq_g} and to the natural logarithm gives
\begin{subequations}
\begin{align}
\nabla[g](x) &=  - \; \sum_{k = 2}^\infty \; \bigg( \sum_{r = 1}^{k-1} \tbinom{k-1}{r} \cdot a_{k-r} \bigg) \cdot x^{-k}    & x > 1/C + 1  \label{eq_Ng_137_4C}\\[8pt]
\nabla[\ln](x) &=\; \Li_1(1/x) & x  > 1  \label{eq_Nln_131A}
\end{align}
\end{subequations}
In \eqref{eq_Nln_131A}, the notation $\Li_s$ refers to the \textit{polylogarithm} \cite{Polylogarithm}.
 \newline\noindent
Finally, declare the operators
\begin{subequations}
\begin{align}
S &:= I - \nabla
\label{eqS} \\[8pt]
J &:= H \circ S
\label{eqJ}
\end{align}
\end{subequations}
In \eqref{eqS}, $I$ denotes the  \textit{identity} operator.

\subsection{3. Decomposition of $f_n$ in terms of $J^{n-k}[1]$ }
It is shown in this section, that solutions $f_n$ of the recurrence relation \eqref{eq_f01_repeat}/\eqref{eq_fn0_repeat}/\eqref{eq_Vermutung_237A}
are linear combinations of $J^{n-k}[1]$.
Applying $\ x \cdot \tfrac{d}{dx}\ $ to $J^m[1]$, one gets
\begin{align}
x \cdot \tfrac{d}{dx} \; J^m[1]\,(x) \;=\; x \cdot \tfrac{d}{dx} \; H\Big[\, S \big[\, &J^{m-1}[1]\,(x) \big] \Big]
\nonumber\\[8pt]
= \hspace{52pt}  S \big[\, &J^{m-1}[1]\,(x) \big]
\nonumber\\[8pt]
= \hspace{70pt} &J^{m-1}[1]\,(x-1)   & x \ge m  \hspace{29pt}  m \ge 1
 \label{eq_dJm}
\end{align}

\vspace{5mm}
\noindent
Hence, replacing $x \rightarrow y+n$ in \eqref{eq_dJm}, for $n \ge m$:
\begin{align}
    (y+n) \; \cdot \; \tfrac{d}{dy} \, J^m[1](y+n) \;&=\; J^{m-1}[1](y+n-1)
    & y \ge 0  \hspace{29pt}  n \ge m \ge 1
 \label{eq_dJm2}
\end{align}

\vspace{1cm}
\noindent
Therefore, if one writes
\begin{align}
f_n(y) &= \sum_{k=0}^n \;\; \beta_k \;\cdot\; J^{n-k}[1]\,(y + n)   & y \ge 0  \hspace{29pt}  n \ge 0
\label{eq_A149}
\end{align}
and applies $(y+n) \cdot \tfrac{d}{dy}$ to \eqref{eq_A149}, while making use of \eqref{eq_dJm2}, one gets
\begin{align}
(y+n) \cdot \tfrac{d}{dy} \; f_n(y)
\;&=\; \sum_{k=0}^n \;\; \beta_k \;\cdot\; (y+n) \cdot \tfrac{d}{dy} \; J^{n-k}[1]\,(y + n)
\nonumber\\[8pt]
  &=\;  \sum_{k=0}^{n-1} \;\; \beta_k \;\cdot\; J^{n-1-k}[1]\,(y+n-1)
\nonumber\\[8pt]
&=\; f_{n-1}(y)   & y \ge 0  \hspace{29pt}  n \ge 1
 \label{eq_dfn}
\end{align}
Here, the fact has been exploited, that $\tfrac{d}{dy} \; J^o[1] \,=\, \tfrac{d}{dy} \;1 \,=\, 0$.
\newline
\eqref{eq_dfn} clearly reproduces \eqref{eq_Vermutung_237A}, while
\eqref{eq_f01_repeat}/\eqref{eq_fn0_repeat} impose the following recurrence relation on the coefficients: \begin{subequations}
\begin{align}
\beta_o &= 1
\label{eq_beta0}\\[8pt]
\beta_n &= -\; \sum_{k=0}^{n-1} \, \beta_k \cdot J^{n-k}[1]\,(n)   &  n \ge 1
\label{eq_betaRec_150A}
\end{align}
\end{subequations}
The first few instances of $\beta_n$ are listed here:
\begin{fleqn} 
\begin{equation}\label{eq_betaList_149}
\footnotesize
\begin{alignedat}{5} 
&\beta_o\; = \;  &&1 \\[3pt]
&\beta_1\; = \;  &&0 \\[3pt]
&\beta_2\; = \; -\; &&J^2[1]\,(2) \\[3pt]
&\beta_3\; = \; -\; &&J^3[1]\,(3) \;+\; J^2[1]\,(2) \cdot J[1]\,(3) \\[3pt]
&\beta_4\; = \; -\; &&J^4[1]\,(4) \;+\; J^2[1]\,(2) \cdot J^2[1]\,(4)
 \; +\; \Big( J^3[1]\,(3) \;-\; J^2[1]\,(2) \cdot J[1]\,(3) \Big) \cdot J[1]\,(4)
\end{alignedat}
\end{equation}
\end{fleqn}
Here, use has been made of the property $J[1]\,(1) = \ln(1) = 0$.

\subsection{4. Decomposition of $J^n[1]$}
Results in this section are stated without giving detailed proofs.
\newline
$J^n[1]$, $n$-fold powers of the $J$ operator \eqref{eqJ} applied to the $1$-function,
may be expanded in terms of powers of the natural logarithm:
\begin{align} \label{eqJm_154B}
J^n[1](x) &= \sum_{j = 0}^n \, (-1)^j \,\cdot\, Q_j(x) \,\cdot\, \frac{\ln(x)^{n-j}}{(n-j)!}    & x \ge n  \hspace{29pt} n \ge 0
\end{align}
 with coefficient functions $Q_j$ fulfilling the recurrence relation
\begin{subequations}
\begin{align}
Q_o \;&=\; 1  \label{eq_Qini}\\[8pt]
Q_{n+1}(x) \;&=\; - H\big[\; \tfrac{1}{n!} \cdot \Li_1(1/x)^n - \nabla[Q_n](x) \,\big]    & x \ge n+1 \hspace{29pt}  n \ge 0 \label{eq_Qrec_160A}
\end{align}
\end{subequations}
Operators $H$ and $\nabla$ are declared in \eqref{eq_H} and \eqref{eq_N} respectively. The notation $\Li_s$ refers to the \textit{polylogarithm} \cite{Polylogarithm}.

\vspace{1cm}
\noindent
The first instances are listed here explicitly:
\begin{subequations}
\begin{align}
Q_o(x) &\;=\; 1    \label{eq_q0_155}\\[8pt]
Q_1(x) &\;=\; 0    \label{eq_q1_155}\\[8pt]
Q_2(x) &\;=\; - H\big[\, \Li_1(1/x) \,\big] \;=\; \Li_2(1/x)   & x \ge 2  \label{eq_q2_155}\\[8pt]
Q_3(x) &\;=\; - H\Big[\; \tfrac{1}{2} \cdot \Li_1(1/x)^2 - \nabla\big[\Li_2(1/x)\big] \,\Big]  \nonumber\\[8pt]
&\;=\;\;\; \sum_{k = 2}^\infty \; \tfrac{1}{k} \cdot \Big( \tfrac{1}{k!}\, \tstirl{k}{2} \,+\, \sum_{r = 1}^{k-1} \tbinom{k-1}{r} \, \tfrac{1}{(k-r)^2} \Big) \cdot x^{-k}
& x \ge 3  \label{eq_q3_155}
\end{align}
\end{subequations}
In \eqref{eq_q3_155}, use has been made of \eqref{eq_Ng_137_4C} and the identity \cite{Flajolet}
\begin{align}\label{eq_flajolet_144A}
    \tfrac{1}{n!} \cdot \Li_1(z)^n &\;=\; \sum_{k = n}^\infty \,  \tfrac{1}{k!} \, \tstirl{k}{n} \cdot z^k
    &  |z| < 1 \hspace{29pt}  n \ge 0 
\end{align}
Here, the notation $\tstirl{k}{n}$ refers to Stirling numbers of the $1^{st}$ kind \cite{Stirling1}.
\newline
For $n \ge 2$, the functions $Q_n(x)$ are power series in $1/x$:
\begin{align} \label{eq_Q}
Q_n(x) &=\; \sum_{s = o}^\infty \, q_{\,n,s} \,\cdot\, x^{-s}    & x \ge n \ge 2
\end{align}
The first non-zero coefficient of $Q_n(x)$ is $q_{n,n-1}$.
In particular, it turns out that
\begin{subequations}
\begin{align}
q_{n+1,\,s} &\;=\; 0  & 0 \le s < n  \hspace{29pt}  n \ge 1
\label{eq_q251B_a}
\\[8pt]
q_{n+1,\,s} &\;=\; \tfrac{1}{s} \cdot \tfrac{1}{s!} \cdot \sum_{\nu = 0}^{n-1}
\;\sum_{\sigma = \nu}^{\nu + s - n} \; \tstirl{s-\sigma}{n-\nu} \cdot \tbinom{s}{\sigma} \cdot A^s_{\sigma,\,\nu}  & s \ge n  \hspace{29pt}  n \ge 1
\label{eq_q251B_b}
\end{align}
\end{subequations}
\eqref{eq_q251B_b} involves Stirling numbers $1^{st}$, binomials, and the quantities $A^s_{m,\,j}$, which are \textit{lower triangular integer} matrices of dimension $s+1$, whose elements obey the recurrence relation
\begin{subequations}
\begin{align}
A^s_{m,\,o} &\;=\; \delta_{m,\,o}  & s \ge m \ge 0
\label{eq_q236A_a}
\\[8pt]
A^s_{m,\,j} &\;=\; \sum_{\mu = 0}^{m-j} \; A^s_{m-\mu-1,\,j-1} \cdot\, \tbinom{m}{\mu+1} \cdot (s-m+1)^{\overline{\,\mu}}  & s \ge m \ge j \ge 1
\label{eq_q236A_b}
\end{align}
\end{subequations}
In \eqref{eq_q236A_a}, $\delta_{i,\,j}$ denotes the Kronecker delta.
\eqref{eq_q236A_b} involves binomials and \textit{rising} factorials.
Special values of $A^s_{m,\,j}$ are
\begin{subequations}
\begin{align}
    A^s_{m,o} \;&=\; \delta_{m,o}  & s \ge m \ge 0
\label{eq_Am0}
\\[8pt]
    A^s_{m,1} \;&=\; (s-1)^{\,\underline{m-1}}  & s \ge m \ge 1
\label{eq_Am1}
\\[8pt]
    A^s_{m,m} \;&=\; m!  & s \ge m \ge 0
\label{eq_Amm}
\end{align}
\end{subequations}
The notation in \eqref{eq_Am1} refers to \textit{falling} factorials.
From the triangle shape of $A^s$ and \eqref{eq_Amm} it follows \cite{A000178}:
\begin{align}
    \det\big(A^s\big) \;&=\; \prod_{k=0}^s k!  & s \ge 0
\label{eq_detA}
\end{align}
Integers in the last row of $A^s$
\begin{ceqn}
\begin{equation}
   A^s_{s,j}  \hspace{30pt} 0 \le j \le s
\label{eq_Asj}
\end{equation}
\end{ceqn}
are related to the Bell matrix with generator $1/j$ for $j>=1$ \cite{BellMatrix} and match \cite{A265607}.
\newline
Example data are listed in Table \ref{list:A}.


\newpage

\begin{table}[ht!]
\footnotesize
\begin{tabular}{r r | r r r r r r r} 
  & & $j$\\
  & & 0 & 1 & 2 & 3 & 4 & 5 & 6\\
 \hline
 $m$ & 0 & 1 \\[8pt]
& 0 & 1 & 0 \\[4pt]
& 1 & 0 & 1 \\[10pt]
& 0 & 1 & 0 & 0 \\[4pt]
& 1 & 0 & 1 & 0 \\[4pt]
& 2 & 0 & 1 & 2 \\[10pt]
& 0 & 1 & 0 & 0 & 0 \\[4pt]
& 1 & 0 & 1 & 0 & 0 \\[4pt]
& 2 & 0 & 2 & 2 & 0 \\[4pt]
& 3 & 0 & 2 & 9 & 6 \\[10pt]
& 0 & 1 & 0 &  0 &  0 &  0 \\[4pt]
& 1 & 0 & 1 &  0 &  0 &  0 \\[4pt]
& 2 & 0 & 3 &  2 &  0 &  0 \\[4pt]
& 3 & 0 & 6 & 15 &  6 &  0 \\[4pt]
& 4 & 0 & 6 & 50 & 72 & 24 \\[10pt]
& 0 & 1 &  0 &   0 &   0 &   0 &   0 \\[4pt]
& 1 & 0 &  1 &   0 &   0 &   0 &   0 \\[4pt]
& 2 & 0 &  4 &   2 &   0 &   0 &   0 \\[4pt]
& 3 & 0 & 12 &  21 &   6 &   0 &   0 \\[4pt]
& 4 & 0 & 24 & 120 & 108 &  24 &   0 \\[4pt]
& 5 & 0 & 24 & 350 & 850 & 600 & 120 \\[10pt]
& 0 & 1 &   0 &    0 &     0 &     0 &    0 &    0 \\[4pt]
& 1 & 0 &   1 &    0 &     0 &     0 &    0 &    0 \\[4pt]
& 2 & 0 &   5 &    2 &     0 &     0 &    0 &    0 \\[4pt]
& 3 & 0 &  20 &   27 &     6 &     0 &    0 &    0 \\[4pt]
& 4 & 0 &  60 &  218 &   144 &    24 &    0 &    0 \\[4pt]
& 5 & 0 & 120 & 1120 &  1750 &   840 &  120 &    0 \\[4pt]
& 6 & 0 & 120 & 3014 & 11250 & 12900 & 5400 &  720 \\[10pt]
& $\vdots$ & & & & $\vdots$
\end{tabular}
\caption{Matrices $A^o$ through $A^6$ computed from \eqref{eq_q236A_a}/\eqref{eq_q236A_b}}
\label{list:A}
\end{table}

\pagebreak




\printbibliography

\end{document}